\documentclass[12pt]{article}
\usepackage{latexsym,epsfig}
\usepackage{amssymb,amsmath,amsthm}

\newcommand{\al}{\alpha}
\newcommand{\be}{\beta}
\newcommand{\ga}{\gamma}

\newcommand{\la}{\lambda}

\renewcommand{\th}{\theta}

\newcommand{\ep}{\epsilon}

\newcommand{\ct}{{CT_p}}
\newcommand{\wct}{{\widetilde{CT}_p}}

\newcommand{\wpi}{{\widehat{\pi}}}

\title{A further note on the inverse nodal problem and Ambarzumyan problem for the $p$-Laplacian}
\author{Y.H.\ Cheng$^1$,\ C.K.\ Law$^2$,\ Wei-Cheng
Lian$^3$ and Wei-Chuan Wang$^4$}
\begin{document}
\maketitle
\begin{abstract}
In this note, we extend some results in a previous paper on the
inverse nodal problem and Ambarzumyan problem for the $p$-Laplacian
to periodic or anti-periodic boundary conditions, and to $L^1$
potentials.
\end{abstract}

%\footnote{Running head: Inverse problems for periodic $p$-Laplacian}

\footnote{$^1$Department of Applied Mathematics, National Sun
Yat-sen University, Kaohsiung, Taiwan 804, R.O.C. E-mail:
jengyh@math.nsysu.edu.tw}
\footnote{$^2$Department of Applied Mathematics, National Sun
Yat-sen University, Kaohsiung, Taiwan 804, R.O.C. E-mail:
law@math.nsysu.edu.tw}
\footnote{$^3$Department of Information Management, National
 Kaohsiung Marine University, Kaohsiung, Taiwan 811, R.O.C.
 E-mail: wclian@mail.nkmu.edu.tw}
\footnote{$^4$Department of Applied Mathematics, National Sun
Yat-sen University, Kaohsiung, Taiwan 804, R.O.C. E-mail:
wangwc@math.nsysu.edu.tw}

\newpage

\def\theequation{\arabic{section}.\arabic{equation}}

\section{Introduction}
\setcounter{equation}{0}
 \hskip0.25in
 Recently, we studied the
$p$-Laplacian with $C^{1}$-potentials and solved the inverse nodal
problem and Ambarzumyan problem for Dirichlet boundary conditions
\cite{LLW08}.
 In this note, we want to extend the results to periodic or anti-periodic boundary conditions, and
to $L^1$ potentials.
 \par
 Consider the equation
\begin{equation}
-\left(y'^{(p-1)}\right)'=(p-1)(\la -q(x))y^{(p-1)}\ ,\label{eq1.1}
\end{equation}
 where $f^{(p-1)}=|f|^{p-1} {\rm sgn} f$. Assume that $q(1+x)=q(x)$ for
 $x\in\mathbb{R}$, then (\ref{eq1.1}) can be coupled with periodic
 or anti-periodic boundary conditions respectively:
\begin{equation}
y(0)=y(1)\ ,\quad y'(0)=y'(1)\qquad \label{eq1.2}\ \ \
\end{equation}
or
\begin{equation}
 y(0)=-y(1)\ ,\quad y'(0)=-y'(1).\qquad
\label{eq1.3}
\end{equation}
When $p=2$, the above is the classical Hill's equation. It follows
from Floquet theory that there are countably many interlacing
periodic and anti-periodic eigenvalues of Hill's operator. However,
Floquet theory does not work for the case $p\neq 2$. In 2001, Zhang
\cite{Zhang2001} studied the properties of eigenvalues for $p>1$
with $L^{1}$-potentials. He applied the rotation number function to
define the minimal eigenvalue $\underline{\lambda}_n(q)$ and the
maximal eigenvalue $\overline{\lambda}_n(q)$  corresponding to
eigenfunctions having $n$ zeros in $[0,1)$, respectively. These
numbers $\underline{\lambda}_n(q)$ and $\overline{\lambda}_n(q)$ are
called rotational periodic eigenvalues and satisfy
\begin{enumerate}
\item[(i)] If $n\in \mathbb{N}\cup \{0\}$ is even, then
$\underline{\lambda}_n(q)$ and $\overline{\lambda}_n(q)$ are
eigenvalues of (\ref{eq1.1}) and (\ref{eq1.2}); if $n\in
\mathbb{N}$ is odd, then $\underline{\lambda}_n(q)$ and
$\overline{\lambda}_n(q)$ are eigenvalues of (\ref{eq1.1}) and
(\ref{eq1.3}).

\item[(ii)]
 $\overline{\lambda}_0(q)<\underline{\lambda}_1(q)\leq\overline{\lambda}_1(q)<\underline{\lambda}_2(q)
\leq \overline{\lambda}_2(q)< \cdots \cdots$.
\end{enumerate}

Although the above properties are very similar to the linear case,
it should be mentioned that the case for the $p$-Laplacian is much
more complicated. For example, for the periodic or anti-periodic
boundary conditions, there  may exist an infinite sequence of
variational eigenvalues and non-variational eigenvalues
(\cite{BR2008}). In the same paper, the authors also showed that the
minimal periodic eigenvalue is simple and variational, while the
minimal anti-periodic eigenvalue is variational but may be not
simple.
 \par
 In 2008, Brown and
Eastham \cite{BE08} derived a sharp asymptotic expansion of
eigenvalues of the $p$-Laplacian with locally integrable and
absolutely continuous $(r-1)$ derivative potentials respectively.
Below is a version of their theorem for periodic eigenvalues of the
$p$-Laplacian (\ref{eq1.1}), (\ref{eq1.2}).

\newtheorem{th1.1}{Theorem}[section]
\begin{th1.1}
\label{th1.1}(\cite[Theorem 3.1]{BE08})
  Let $q$ be $1$-periodic and
locally integrable in $(-\infty, \infty)$.
 %Denote by $\gamma_{n}$ and $c_{n}$  the complex
 %Fourier coefficients of $|S_p(x)|^{p}-\frac{1}{p}$ and $q$ on
 %$(0,1)$ respectively.
 Then the rotationally periodic eigenvalue
$\la_{2n}=\underline{\lambda}_{2n}$, or $\overline{\lambda}_{2n}$
satisfies
\begin{equation}
\label{eq1.4}
\lambda_{2n}^{1/p}=2n\widehat{\pi}+\frac{1}{p(2n\widehat{\pi})^{p-1}}
\int_{0}^{1}q(t)dt+o(\frac{1}{n^{p-1}}).
 %\frac{1}{(p-1)(2n\widehat{\pi})^{p-1}}\delta_{n}+O(\frac{1}{n^{2(p-1)}})\
 %,
\end{equation}
 %where $|\delta_{n}|\leq \sum_{k=-\infty}^{\infty}
 %|\gamma_{k}||c_{2k(n+1)}|=o(1)$ as $n\to\infty$.
\end{th1.1}

By a similar argument, the asymptotic expansion of the anti-periodic
eigenvalue ${\lambda}_{2n-1}=\underline{\lambda}_{2n-1}$ or
 $\overline{\lambda}_{2n-1}$, which corresponds to the anti-periodic
eigenfunction with $2n-1$ zeros in $[0,1)$, satisfies
\begin{equation}
\label{eq1.5}
\lambda_{2n-1}^{1/p}=(2n-1)\widehat{\pi}+\frac{1}{p((2n-1)\widehat{\pi})^{p-1}}
\int_{0}^{1}q(t)dt+o(\frac{1}{n^{p-1}})
 %\frac{1}{(p-1)((2n-1)\widehat{\pi})^{p-1}}\delta_{n}+O(\frac{1}{n^{2(p-1)}})\
.\end{equation}
 \hskip0.25in
 The inverse nodal problem is the problem of understanding
the potential function through its nodal data. In 2006, some of us
(C.-L.) \cite{CL062} studied  Hill's equation. We first made a
translation of the interval by the first nodal length so that the
periodic problem is reduced to a Dirichlet problem, and then solved
the uniqueness, reconstruction and stability problems using the
nodal set of periodic eigenfunctions.
 \par
 We denote by $\{
x_i^{(n)}\}_{i=0}^{n-1}$ the zeros of the eigenfunction
corresponding to $\la_{n}$, and define the nodal length
$\ell_{i}^{(n)}=x_{i+1}^{(n)}-x_{i}^{(n)}$ and $j=j_{n}(x)=\max\{i
:\ x_{i}^{(n)}\leq x \}$. Our main theorem is as follows.

\newtheorem{th1.2}[th1.1]{Theorem}
\begin{th1.2}
\label{th1.2}
 Let $q\in L^{1}(0,1)$ be $1$-periodic. Define  $F_n(x)$
as the following:
\begin{enumerate}
\item[(a)]For  periodic boundary condition, let
$$F_{2n}(x)=p(2n\widehat{\pi})^p[(2n)\ell_j^{(2n)}-1]+\int_0^1q(t)dt,$$
\item[(b)]For the anti-periodic boundary condition, let
$$F_{2n-1}(x)=p((2n-1)\widehat{\pi})^p[(2n-1)\ell_j^{(2n-1)}-1]+\int_0^1q(t)dt.$$
\end{enumerate}

Then both $\{F_{2n}\}$ and $\{F_{2n-1}\}$  converges to $q$
pointwisely a.e. and in $L^{1}(0,1)$.
\end{th1.2}

Thus either one of the sequences $\{F_{2n}\}/ \{F_{2n-1}\}$ will
give the reconstruction formula for $q$. Note that here $q\in
L^1(0,1)$. Furthermore, the map between the nodal space and the
set of admissible potentials are homeomorphic after a partition
(cf.\cite{LLW08}). The same idea
 also works for linear separated boundary value problems with
 integrable potentials.
 \par
 Using the eigenvalue
asymptotics above, the Ambarzumyan problems for the periodic and
anti-periodic boundary conditions can also be solved.
\newtheorem{th1.3}[th1.1]{Theorem}
\begin{th1.3}
\label{th1.3}Let $q\in L^{1}(0,1)$ be periodic of period $1$.
\begin{enumerate}
\item[(a)] If the spectrum of periodic eigenvalues problem (\ref{eq1.1}), (\ref{eq1.2}) contains $\{(2n\widehat{\pi})^p:n\in \mathbb{N}\cup \{0\}\}$
and $0$ is the least eigenvalue, then $q=0$ on $[0,1]$.
\item[(b)] If the spectrum of anti-periodic eigenvalue problem (\ref{eq1.1}), (\ref{eq1.3}) contains $\{((2n-1)\widehat{\pi})^p:n\in
\mathbb{N}\}$;  $\widehat{\pi}^p$ is the least eigenvalue and
$\int_0^1q(t)(S_p(\widehat{\pi}t)S_p'(\widehat{\pi}t)^{(p-1)})'dt=0$,
then $q=0$ on $[0,1]$.
\end{enumerate}
\end{th1.3}

 In section 2, we shall apply Theorem \ref{th1.1} to study
 on periodic and
 anti-periodic boundary conditions. In section 3, we shall deal with
 the case of
 linear separated boundary
 conditions.
 \par
 The stability issue of the inverse nodal problem with $L^1$
 potentials associated with perodic/antiperiodic as well as linear
 separated boundary conditions can also be proved.  The proof goes
 in the same manner as in \cite{LLW08} and is so omitted.
\section{ Proof of main results}
\setcounter{equation}{0}
 \hskip0.25in
 Fix $p>1$ and assume that $q=0$ and $\lambda=1$. Then
(\ref{eq1.1}) becomes
$$-(y'^{(p-1)})'=(p-1)y^{(p-1)}.$$
 Let $S_p$ be the solution
satisfying the initial conditions $\displaystyle S_p(0)=0$,
$S_p'(0)=1$. It is well known that  $S_p$ and its derivative $S_p'$
are periodic functions on $\mathbb{R}$ with period $2\widehat\pi$,
where $ {\widehat \pi}=\frac{2\pi}{p\sin(\frac{\pi}{p})}$. The two
functions also satisfy the following identities (cf.\
\cite{BE08,LLW08}).
\newtheorem{th2.1}{Lemma}[section]
\begin{th2.1}
\label{th2.1}
\begin{enumerate}
\item[(a)] $|S_p(x)|^p+|S_p'(x)|^p=1$ for any $x\in \mathbb{R}$;
\item[(b)] $(S_pS_p'^{(p-1)})'=|S_p'|^p-(p-1)|S_p|^p=1-p|S_p|^p=(1-p)+p|S_p'|^p$\ .
\end{enumerate}
\end{th2.1}
 Next we define a generalized Pr\"ufer substitution using $S_p$ and
 $S_p'$:
\begin{equation}
y(x)=r (x)S_p(\lambda ^{1/p}\theta (x)), \ \ y'(x)=\lambda ^{1/p}
r(x)S_p'(\lambda ^{1/p}\theta (x))\ . \label{eq2.1}
\end{equation}
 By Lemma
\ref{th2.1}, one obtains (\cite{LLW08})
\begin{equation}
\theta '(x)=1-\frac {q}{\lambda}|S_p(\lambda ^{1/p}\theta (x))|^p\ .
\label{eq2.2}
\end{equation}

\newtheorem{th3.1}[th2.1]{Theorem}
\begin{th3.1}
\label{th3.1}
 In the periodic/antiperiodic eigenvalue problem, if $q\in L^{1}(0,1)$ be periodic of period $1$, then
 $$ q(x)=\lim _{n\rightarrow \infty}p\lambda _n\left(
\frac{\lambda _n^{1/p}\ell_j^{(n)}}{\widehat{\pi}}-1\right)\ ,
 $$ pointwisely  a.e. and in $L^{1}(0,1)$, where $j=j_n(x)=\max\{k:x_k^{(n)}\leq x\}$.
\end{th3.1}

 The proof below works for both even and odd $n$'s, i.e.\ for both periodic and antiperiodic problems.
  Some of the arguments above are motivated by
\cite{CW07}. See also \cite{LSY1999}.
\begin{proof}
  First, integrating (\ref{eq2.2}) from $x_{k}^{(n)}$
to $x_{k+1}^{(n)}$ with $\la=\la_n$, we have
\begin{eqnarray*}\frac{\widehat{\pi}}{\lambda_{n}^{1/p}}&=&\ell_{k}^{(n)}-\int_{x_{k}^{(n)}}^{x_{k+1}^{(n)}}\frac{q(t)}{\la_{n}}|S_p(\lambda_{n}^{1/p}\theta(t))|^{p}dt\
,\\
&=&\ell_{k}^{(n)}-\frac{1}{p\la_{n}}\int_{x_{k}^{(n)}}^{x_{k+1}^{(n)}}q(t)dt-\frac{1}{\la_{n}}\int_{x_{k}^{(n)}}^{x_{k+1}^{(n)}}q(t)(|S_p(\lambda_{n}^{1/p}\theta(t))|^{p}-\frac{1}{p})dt\
.
\end{eqnarray*}
Hence,
\begin{equation}
\ell_{k}^{(n)}=\frac{\widehat{\pi}}{\lambda_{n}^{1/p}}
+\frac{1}{p\la_{n}}\int_{x_{k}^{(n)}}^{x_{k+1}^{(n)}}q(t)dt+\frac{1}{\la_{n}}\int_{x_{k}^{(n)}}^{
x_{k+1}^{(n)}}q(t)(|S_p(\lambda_{n}^{1/p}\theta(t))|^{p}-\frac{1}{p})dt\
.
\end{equation}
and
\begin{equation}
 p\lambda_n\left(
\frac{\lambda _n^{1/p}\ell_k^{(n)}}{\widehat{\pi}}-1\right)=\frac{\lambda _n^{1/p}}{
\widehat{\pi}}\int_{x_{k}^{(n)}}^{x_{k+1}^{(n)}}q(t)dt+\frac{p\lambda _n^{1/p}}{\widehat{\pi}}
\int_{x_{k}^{(n)}}^{x_{k+1}^{(n)}}q(t)(|S_p(\lambda_{n}^{1/p}\theta(t))|^{p}-\frac{1}{p})dt\
\label{eq2.3} .\end{equation}

Now, for $x\in (0,1)$, let $j=j_{n}(x)=\max\{k : x_{k}^{(n)}\leq x
\}$. Then $x\in [x_{j}^{(n)},x_{j+1}^{(n)})$ and, for large $n$,
$$[x_{j}^{(n)},x_{j+1}^{(n)}) \subset
B(x,\frac{2\widehat{\pi}}{\lambda_{n}^{1/p}})\ ,$$ where
$B(t,\varepsilon)$ is the open ball centering $t$ with radius
$\varepsilon$. That is, the sequence of intervals
$\{[x_{j}^{(n)},x_{j+1}^{(n)}) : n
\mbox{ is sufficiently large}\}$ shrinks  to $x$ nicely (cf. Rudin
\cite[p.140]{Rudin}). Since $q\in L^{1}(0,1)$ and
$\frac{\lambda _n^{1/p}\ell_{k}^{(n)}}{\widehat{\pi}}=1+o(1)$, we
have $$h_{n}(x)\equiv\frac{\lambda
_n^{1/p}}{\widehat{\pi}}\int_{x_{j}^{(n)}}^{x_{j+1}^{(n)}}q(t)dt =
\frac{\lambda _n^{1/p}\ell_{j}^{(n)}}{\widehat{\pi}}
\frac{1}{\ell_{j}^{(n)}}\int_{x_{j}^{(n)}}^{x_{j+1}^{(n)}}q(t)dt $$
converges to $q(x)$ pointwisely a.e. $x\in (0,1)$. Furthermore,
since $$|h_{n}(x)|\leq \frac{\lambda
_n^{1/p}}{\widehat{\pi}}\int_{x_{j}^{(n)}}^{x_{j+1}^{(n)}}|q(t)|dt\equiv
g_{n}(x)\ ,$$ and
$$\int_{0}^{1}g_{n}(t)dt=\sum_{k=0}^{n-1}\frac{\lambda
_n^{1/p}\ell_{k}^{(n)}}{\widehat{\pi}}\int_{x_{k}^{(n)}}^{x_{k+1}^{(n)}}|q(t)|dt=(1+o(1))\|q\|_{1}\
, $$ we have $ h_{n}(t)\to q(t)$ in $L^{1}(0,1)$ by Lebesgue
dominated convergence theorem.
 On the other hand, let
$q_{k,n}\equiv\frac{1}{\ell_{k}^{(n)}}\int_{x_{k}^{(n)}}^{x_{k+1}^{(n)}}q(t)dt$.
Then $q_{j,n}$ converges to $q$ pointwisely a.e. $x\in (0,1)$. Let
$\phi_{n}(t)=|S_p(\lambda_{n}^{1/p}\theta(t))|^{p}-\frac{1}{p}$.
 Then
\begin{eqnarray*}
T_{n}(x)&\equiv&\frac{p\lambda
_n^{1/p}}{\widehat{\pi}}\int_{x_{j}^{(n)}}^{x_{j+1}^{(n)}}q(t)\phi_{n}(t)
dt\
,\\
&=&\frac{p\lambda
_n^{1/p}}{\widehat{\pi}}\int_{x_{j}^{(n)}}^{x_{j+1}^{(n)}}(q(t)-q_{j,n})\phi_{n}(t)dt
+\frac{p\lambda
_n^{1/p}}{\widehat{\pi}}\int_{x_{j}^{(n)}}^{x_{j+1}^{(n)}}q_{j,n}\phi_{n}(t)dt\
,\\
&\equiv&A_{n}+B_{n}\ .
\end{eqnarray*}
By Lemma \ref{th2.1}(b)
 and  (\ref{eq2.2}),
\begin{eqnarray*}B_{n}
 % &=&\frac{p\lambda
 %_n^{1/p}}{\widehat{\pi}}q_{j,n}\int_{x_{j}^{(n)}}^{x_{j+1}^{(n)}}U(t)dt\
 %,\\
&=&\frac{p\lambda
_n^{1/p}q_{j,n}}{\widehat{\pi}}\int_{x_{j}^{(n)}}^{x_{j+1}^{(n)}}\left(
|S_p(\lambda_{n}^{1/p}\theta(t))|^{p}-\frac{1}{p}
\right)\left(\theta'(t)+
\frac{q(t)}{\lambda_{n}}|S_p(\lambda_{n}^{1/p}\theta (t))|^p\,
\right)\, dt
,\\
&=&-\left.\frac{p q_{j,n}}{\widehat{\pi}}S_p(\la_n^{1/p}\th(t))
S_p'(\la_{n}^{1/p}\th(t))^{(p-1)}\right|_{x_{j}^{(n)}}^{x_{j+1}^{(n)}}
+O(\la_{n}^{-1+1/p})\ ,\\
&=&O(\la_{n}^{-1+1/p})\ .
\end{eqnarray*}
Also,
 \begin{eqnarray*}|A_{n}|&\leq&
\frac{p\lambda_n^{1/p}}{\widehat{\pi}}\int_{x_{j}^{(n)}}^{x_{j+1}^{(n)}}
|q(t)-q_{j,n}|||S_p(\lambda_{n}^{1/p}\theta(t))|^{p}-\frac{1}{p}|dt\ ,\\
&\leq&
\frac{(p-1)\lambda_n^{1/p}}{\widehat{\pi}}\int_{x_{j}^{(n)}}^{x_{j+1}^{(n)}}|q(t)-q_{j,n}|dt
\ ,
\end{eqnarray*}
which converges to $0$ pointwisely a.e.\ $x\in (0,1)$ because the
sequence of intervals $\{[x_{j}^{(n)},x_{j+1}^{(n)}) : n
\mbox{ is sufficiently large}\}$ shrinks  to $x$ nicely. We
conclude that $T_{n}(x)\to 0$ a.e. $x\in (0,1)$. Finally,
applying Lebesgue dominated convergence theorem as above,
$T_{n}(x)\to 0$ in $L^{1}(0,1)$. Hence the left hand side of
(\ref{eq2.3}) converges to $q$ pointwisely  a.e. and in
$L^{1}(0,1)$.
\end{proof}
 \vskip0.1in
\begin{proof}[Proof of Theorem \ref{th1.2}]
 \hfill

By the eigenvalue
estimates (\ref{eq1.4})   and (\ref{eq1.5}), we have
\begin{equation}
p\lambda_{2n}(\frac{\lambda_{2n}^{1/p}\ell_{j_{2n}(x)}^{(2n)}}{\widehat{\pi}}-1)
=p(2n\widehat{\pi})^p(2n\ell_j^{(2n)}-1)+2n\ell_{j_{2n}(x)}^{(2n)}\int_0^1q(t)dt+o(1)\
.\label{eq3.3}
\end{equation}
 Hence by Theorem \ref{th3.1} and the fact that
$2n\ell_j^{(2n)}=1+o(1)$,
 $$
 F_{2n}(x)\equiv
 p(2n\widehat{\pi})^p(2n\ell_j^{(2n)}-1)+\int_0^1q(t)dt$$
 also converges
to $q$ pointwisely a.e. and in $L^{1}(0,1)$. The proof for (b) is
the same.
\end{proof}

\begin{proof}[Proof of Theorem \ref{th1.3}]
 \hfill

Here we only give the proof of (b). First, since all anti-periodic
eigenvalues include $\{((2n-1)\widehat{\pi})^p:n\in \mathbb{N}\}$,
we have, by (\ref{eq1.5}), $\int_0^1q(t)dt=0$.

Moreover,  $S_p(\widehat{\pi}x)$ satisfies anti-periodic boundary
conditions. So by Lemma \ref{th2.1}(b),
 $$
\int_0^1|S_p'(\widehat{\pi}t)|^pdt-\frac{p-1}{p} =
\int_0^1q(t)|S_p(\widehat{\pi}t)|^pdt=
\int_0^1|S_p(\widehat{\pi}t)|^pdt-\frac{1}{p}=0\ .
 $$
Hence, by  the variational principle, we have
 $$
\widehat{\pi}^p=\lambda_1
\leq
\frac{\int_0^1\widehat{\pi}^p|S_p'(\widehat{\pi}t)|^pdt+(p-1)\int_0^1q(t)|S_p(\widehat{\pi}t)|^pdt}
{(p-1)\int_0^1|S_p(\widehat{\pi}t)|^pdt}= \widehat{\pi}^p\ .
 $$
 This
implies $S_p(\widehat{\pi}x)$ is the first eigenfunction. Therefore
$q=0$
 on $[0,1]$.
\end{proof}

\section{Linear separated boundary conditions}
\setcounter{equation}{0}
\hskip0.25in
 Consider the one-dimensional $p$-Laplacian with linear separated
 boundary conditions
  \begin{equation}
  \left\{ \begin{array}{l}
   y(0)S_p'(\al)+y'(0)S_p(\al) = 0\\
   y(1) S_p'(\be)+y'(1)S_p(\be)= 0
   \end{array}
   \right.\ ,\label{eq4.1}
  \end{equation}
  where $\al,\be\in [0,\widehat{\pi})$. Letting $\la_n$ be the $n$th eigenvalue whose associated eigenfunction has
  exactly $n-1$ zeros in $(0,1)$, the generalized phase $\th_n$ as given in (\ref{eq2.2})
  satisfies
  \begin{equation}
  \th_n(0)=\frac{-1}{\la_n^{1/p}}\wct^{-1}(-\frac{\wct(\al)}{\la_n^{1/p}});\qquad
  \th_n(1)=\frac{1}{\la_n^{1/p}}\left( n\widehat{\pi}-
  \wct^{-1}(-\frac{\wct(\be)}{\la_n^{1/p}})\right)\ ,\label{eq4.8}
  \end{equation}
  where the function $\ct(\ga):=\frac{S_p(\ga)}{S_p'(\ga)}$ is an analogue of
  cotangent function, while $\wct(\ga):=\ct(\ga)$ if $\ga\neq 0$; and
  $\wct(\ga):=0$ otherwise. Also $\wct^{-1}$ stands for the inverse
  of $\wct$, taking values only in $[0,\widehat{\pi})$.
  \par
 Let $\phi_{n}(x)=|S_p(\lambda_{n}^{1/p}\th_n(x)|^p-\frac{1}{p}$,
  where .
  Below we shall state a  general Riemann-Lebesgue lemma,
 which shows that $\int_0^1\!\phi_n g\rightarrow 0$ for any $g\in L^1(0,1)$, when $\la_n$'s are associated
 with a certain linear separated boundary conditions.
 In the case of periodic boundary conditions,
 Brown and Eastham \cite{BE08} used a Fourier series expansion of $\phi_{n}$ where
 $\phi_{n}(\la_n^{1/p}\th_n(x))\approx \phi_{n}(\al+2n\widehat{\pi} x)$ and apply
 Plancherel Theorem to show convergence.

 \newtheorem{th2.2}{Lemma}[section]
 \begin{th2.2}
 \label{th2.2}
 Let $f_{n}$ be uniformly bounded and integrable on $(0,1)$.
 Suppose for each $n$, there exists a partition $\{
 x_0^n=0<x_1^n<\cdots <x_n^n=1\}$ such that $\Delta
 x_k^n=o(1)$, and $F_k^n(x):=\int_{x_k^n}^x f_{n}(t)\, dt$
 satisfies $F_k^n(x)=O(\frac{1}{n})$ for $x\in (x_{k}^{n},x_{k+1}^{n})$ and
 $F_k^n(x_{k+1}^n)=o(\frac{1}{n})$ uniformly in $k=1,\ldots,n-2$, as
 $n\to\infty$. Then for any $g\in L^1(0,1)$, $\int_0^1
 gf_{n}\rightarrow 0$ as $n\to\infty$.
 \end{th2.2}
 \begin{proof}
 Take any $\ep>0$, there is a $C^1$ function $\tilde{g}$ on $[0,1]$
 such that $\int_0^1|\tilde{g}-g|<\ep$.  Let $|f_{n}|, |\tilde{g}|\leq M$.
 Then
 $$
 \int_0^1 g f_{n} =\int_0^1 (g-\tilde{g})f_{n} +\int_0^1 \tilde{g} f_{n},
 $$
 where $|\int_0^1 (g-\tilde{g})f_{n}|\leq M\ep$. Also
 $$
 \int_0^1 \tilde{g}f_{n}=\sum_{k=0}^{n-1}\int_{x_k^n}^{x_{k+1}^n}
 \tilde{g}f_{n}=\sum_{k=1}^{n-2}
 \left(\tilde{g}(x_{k+1}^n)F(x_{k+1}^n)-\int_{x_k^n}^{x_{k+1}^n}
 \tilde{g}' F_k^n \right)+o(1),
 $$
 where
 $$
 |\int_{x_k^n}^{x_{k+1}^n}
 \tilde{g}' F_k^n|=O(\frac{1}{n})\, \int_{x_k^n}^{x_{k+1}^n}
 |\tilde{g}'|=o(\frac{1}{n}).
 $$
 Therefore $\int_0^1 \tilde{g}f_{n}=o(1)$ as $n\to\infty$.
 \end{proof}
 \newtheorem{th2.3}[th2.2]{Corollary}
 \begin{th2.3}
 \label{th2.3}
 Consider the $p$-Laplacian (\ref{eq1.1}) with boundary conditions
 (\ref{eq4.1}).
 Define $\phi_{n}(x)=|S_p(\lambda_{n}^{1/p}\th_n(x))|^p-\frac{1}{p}$, then
 for any $g\in L^1(0,1)$, $\int_0^1 \phi_{n}g\rightarrow 0$.
 \end{th2.3}
 \begin{proof}
 Since $\th_n(0)$ and $\th_n(1)$ are as given in (\ref{eq4.8}),
   $\phi_{n}$ is
 uniformly bounded on $[0,1]$. Take $x_k^n$ be such that
 $\th(x_k^n)=\frac{k\widehat{\pi}}{\lambda_{n}^{1/p}}$.
 Also by integrating the phase equation (\ref{eq2.2}),
 $\la_{n}^{1/p}=O(n)$, and
  $$
  \Delta x_n=O(\frac{1}{\la_{n}^{1/p}})=O(\frac{1}{n}).
  $$
  Hence by Lemma \ref{th2.1}(b) and (\ref{th2.2}), we have for
  $k=1,\ldots,n-2$,
 \begin{eqnarray*}
 \int_{x_k^n}^{x_{k+1}^n} \phi_{n}(x)\, dx &=& \frac{-1}{p\la_n^{1/p}}\int_{x_k^n}^{x_{k+1}^n}
  \frac{1}{\th'_n(x)}\, \frac{d}{dx}\left[
  S_p(\la_n^{1/p}\th_n(x))S_p'(\la_n^{1/p}\th_n(x))^{(p-1)}\right]\,
  dx\ ,\\
  &=& \frac{-1}{p\la_n^{1/p}} \left[
  S_p(\la_n^{1/p}\th_n(x))S_p'(\la_n^{1/p}\th_n(x))^{(p-1)}\right]_{x_k^n}^{x_{k+1}^n}
  +O(\frac{1}{\la_{n}})\ ,\\
  &=&O(\frac{1}{\la_{n}})=o(\frac{1}{n})\ ,
  \end{eqnarray*}
  since $S_p(k\widehat{\pi})=0$. It is also clear that
  $\int_{x_k^n}^{x} \phi_{n}(x)\, dx =O(\frac{1}{n})$. Thus we may apply
  Lemma~\ref{th2.2} to complete the proof.
  \end{proof}
 \newtheorem{th4.1}[th2.2]{Theorem}
\begin{th4.1}
\label{th4.1}
 When $q\in L^1(0,1)$,
 the eigenvalues $\lambda_n$ of the Dirichlet $p$-Laplacian (\ref{eq1.1})
 satisfies, as $n\rightarrow \infty$,
 \begin{equation}
 \lambda _n^{1/p}
 = n{\widehat \pi} + \frac{1}{p(n{\widehat \pi})^{p-1}} \int_0^1
 q(t)dt+o(\frac{1}{n^{p-1}})\ .\label{eq4.3}
  \end{equation}
Furthermore, $F_n$ converges to $q$ pointwisely and in $L^1(0,1)$,
where
 $$F_n(x):=p(n\widehat{\pi}
 )^p(n\ell_j^{(n)}-1)+\int _0^1 q(t)\, dt.
 $$
 \end{th4.1}
 \begin{proof}
 Integrating (\ref{eq2.2}) from $0$ to $1$, we have
\begin{eqnarray*}
\la_{n}^{1/p}&=&n\widehat{\pi}
+\frac{1}{p\la_{n}^{1-1/p}}\int_{0}^{1}q(t)|S_p(\la_{n}^{1/p}\theta
(t))|^{p}dt\ ,\\
 &=& n\widehat{\pi}+\frac{1}{p\la_{n}^{1-1/p}}\int_{0}^{1}q(t)dt
+\frac{1}{p\la_{n}^{1-1/p}}\int_{0}^{1}q(t)(|S_p(\la_{n}^{1/p}\theta
(t))|^{p}-\frac{1}{p})dt\ .
\end{eqnarray*}
Then by Corollary \ref{th2.3}, we have
 $$
 \int_{0}^{1}q(t)(|S_p(\la_{n}^{1/p}\theta
(t))|^{p}-\frac{1}{p})dt=o(1)\ , $$ for any $q\in L^{1}(0,1)$. Hence
(\ref{eq4.3}) holds. Furthermore, by Theorem \ref{th3.1}, we can
obtain the reconstruction formula with pointwise and $L^{1}$
convergence.
\end{proof}
\noindent {\bf Remark.} In the same way, the Ambarzumyan Theorems
for Neumann as well as Dirichlet boundary conditions as given in
\cite[Theorems 1.3 and 5.1]{LLW08} can also be extended to work for
$L^1$ potentials. On the other hand,  for general linear separated
boundary problems (\ref{eq4.1}),
 \begin{equation}
 \lambda_n^{1/p}=n_{\al\be}\widehat{\pi}+\frac{(\wct(\be))^{(p-1)}-(\wct(\al))^{(p-1)}}{(n_{\al\be}\widehat{\pi})^{p-1}}
 +\frac{1}{p(n_{\al\be}\widehat{\pi})^{p-1}}\int
_0^1q(x)\, dx+o(\frac{1}{n^{p-1}}), \label{eq4.11}
  \end{equation}
  where
 $$
 n_{\al\be}=\left\{ \begin{array}{ll}
  n & \mbox{if }\al=\be=0\\
  n-1/2 & \mbox{if }\al>0=\be \mbox{ or } \be>0=\al \\
  n-1 & \al,\be>0
  \end{array}
  \right.
  $$
  This is because, after an integration of (\ref{eq2.2}),
  \begin{equation}
  \th_n(1)-\th_n(0)=1-\frac{1}{\la_n}\int_0^1
  q(x)|S_p(\la_n^{1/p}\th(x))|^p\,
  dx+o(\frac{1}{\la_n}).\label{eq4.9}
  \end{equation}
  By (\ref{eq4.8}), if $\al=0$, then $\th_n(0)=0$. Similarly
  $\th_n(1)=0$ if $\beta=0$.
Now,  let $y=CT_{p}^{-1}(x)$. Then
$x=CT_{p}(y)$ and hence  
 $$y'=
-\frac{1/x^{2}}{1+\frac{1}{|x|^{p}}}
=\frac{-|x|^{p-2}}{1+|x|^{p}}=-|x|^{p-2}(1+O(|x|^{p}),
 $$ 
 when $|x|$ is sufficiently small. Since $y(0)=\frac{\wpi}{2}$,
we have  $$y(x)=\frac{\wpi}{2} -
\frac{x^{(p-1)}}{p-1}+O(x^{2p-1})\ .$$
   Therefore, when $n$ is sufficiently large,
   $$
   \th_n(0)=\frac{\wpi}{2\la_n^{1/p}}+\frac{(\ct(\al))^{(p-1)}}{(p-1)\la_n^{(p-1)/p}}+O(\la_{n}^{\frac{1-2p}{p}}).
   $$
   Similarly, when $\be\neq 0$,
   $$
   \th_n(1)=\frac{(n-\frac{1}{2})\wpi}{\la_n^{1/p}}+\frac{(\ct(\be))^{(p-1)}}{(p-1)\la_n^{(p-1)/p}}+O(\la_{n}^{\frac{1-2p}{p}}).
   $$
   Hence (\ref{eq4.11}) is valid.
   Furthermore, $F_n$ converges to $q$ pointwisely and in $L^1(0,1)$,
 where
 $$
 F_n(x):=p(n_{\al\be}\widehat{\pi})^p
 \left[ (n_{\al\be}+\frac{(\wct(\be))^{(p-1)}-(\wct(\al))^{(p-1)}}{(n_{\al\be}\widehat{\pi})^{p-1}}
 )\ell_j^{(n)}-1\right]+\int _0^1
 q(t)\, dt.
 $$

\section*{Acknowledgments}
 \hskip0.25in
The authors are supported in part by National Science Council,
Taiwan under contract numbers NSC 98-2115-M-110-006, NSC
97-2115-M-005-MY2 and  NSC 97-2115-M-022-001.

\end{document}